\documentclass[11pt]{article}
\usepackage[matrix,arrow,curve,cmtip]{xy}
\usepackage{amssymb}
\usepackage{theorem}
\usepackage{latexsym}
\usepackage{graphicx}
\usepackage{theorem}
\usepackage[a4paper,width=140mm,height=220mm]{geometry}

\def\to{\mbox{$\xymatrix@1@C=5mm{\ar@{->}[r]&}$}}
\def\tto{\mbox{$\xymatrix@1@C=5mm{\ar@{=>}[r]&}$}}

\def\bkar{\mbox{$\xymatrix@1@C=5mm{\ar@{->}[l]&}$}}
\def\distsign{\begin{picture}(0,0)\put(0,0){\circle{4}}\end{picture}}
\def\dist{\mbox{$\xymatrix@1@C=5mm{\ar@{->}[r]|{\distsign}&}$}}
\def\bkdist{\mbox{$\xymatrix@1@C=5mm{\ar@{->}[l]|{\distsign}&}$}}
\def\biar{\mbox{$\xymatrix@1@C=5mm{\ar@<1.5mm>[r]\ar@<-0.5mm>[r]&}$}}
\def\bidist{\mbox{$\xymatrix@1@C=5mm{\ar@<1.5mm>[r]|{\distsign}\ar@<-0.5mm>[r]|{\distsign}&}$}}
\def\adjar{\mbox{$\xymatrix@1@C=5mm{\ar@<1.5mm>@{<-}[r]\ar@<-0.5mm>[r]&}$}}
\def\adjdist{\mbox{$\xymatrix@1@C=5mm{\ar@<1.5mm>@{<-}[r]|{\distsign}\ar@<-0.5mm>[r]|{\distsign}&}$}}
\def\iso{\mbox{$\xymatrix@1@C=6mm{\ar@{->}[r]^{\sim}&}$}}
\def\doubiso{\mbox{$\xymatrix@1@C=6mm{\ar@{<->}[r]^{\sim}&}$}}
\def\doubar{\mbox{$\xymatrix@1@C=6mm{\ar@{<->}[r]&}$}}

\font\atip = xycmat12 at 10pt
\font\btip = xycmbt12 at 10pt
\def\tip#1{
\mbox{
\begin{picture}(0,0)
\put(-725,-25){\atip\char #1}
\put(-725,-25){\btip\char #1}
\end{picture} }}
\newsavebox{\westsouthwesthead}
\savebox{\westsouthwesthead}{%
\tip{119}}
\newcommand{\wswhead}{\usebox{\westsouthwesthead}}
\def\Endoar#1{
\setlength{\unitlength}{0.01pt}
\ifinner
\mbox{
\begin{picture}(300,1200)(1300,0)
\put(1450,680){\mbox{\footnotesize{$#1$}}}
\put(900,770){\oval(900,900)[t]}
\put(900,770){\oval(900,900)[br]}
\put(900,300){\wswhead}
\end{picture}}
\else
\mbox{
\begin{picture}(200,1400)(1600,400)
\put(2100,1300){\mbox{$#1$}}
\put(1300,1300){\oval(1300,1300)[t]}
\put(1300,1300){\oval(1300,1300)[br]}
\put(1300,600){\wswhead}
\end{picture}}
\fi
\setlength{\unitlength}{1pt}}
\def\endoar{
\setlength{\unitlength}{0.01pt}
\ifinner
\mbox{
\begin{picture}(300,1200)(1300,0)
\put(900,770){\oval(900,900)[t]}
\put(900,770){\oval(900,900)[br]}
\put(900,300){\wswhead}
\end{picture}}
\else
\mbox{
\begin{picture}(200,1400)(1600,400)
\put(1300,1300){\oval(1300,1300)[t]}
\put(1300,1300){\oval(1300,1300)[br]}
\put(1300,600){\wswhead}
\end{picture}}
\fi
\setlength{\unitlength}{1pt}}
\def\Endodist#1{
\setlength{\unitlength}{0.01pt}
\ifinner
\mbox{
\begin{picture}(300,1200)(1300,0)
\put(1600,770){\mbox{\footnotesize{$#1$}}}
\put(900,770){\oval(900,900)[t]}
\put(900,770){\oval(900,900)[br]}
\put(1300,970){\circle{400}}
\put(900,300){\wswhead}
\end{picture}}
\else
\mbox{
\begin{picture}(200,1400)(1600,400)
\put(2100,1300){\mbox{$#1$}}
\put(1300,1300){\oval(1300,1300)[t]}
\put(1300,1300){\oval(1300,1300)[br]}
\put(1850,1600){\circle{400}}
\put(1300,600){\wswhead}
\end{picture}}
\fi
\setlength{\unitlength}{1pt}}

\newtheorem{theorem}{Theorem}[section]
\newtheorem{lemma}[theorem]{Lemma}
\newtheorem{definition}[theorem]{Definition} 
\newtheorem{proposition}[theorem]{Proposition}
\newtheorem{corollary}[theorem]{Corollary}

{\theorembodyfont{\upshape}\newtheorem{example}[theorem]{Example}}
\newcommand{\proof}{\noindent {\it Proof\ }: }
\def\endofproof{$\mbox{ }\hfill\Box$\par\vspace{1.8mm}\noindent}

\def\Psd{{\sf Psd}}
\def\ClPsd{{\sf ClPsd}}
\def\lax{_{\sf lax}}
\def\Tens{{\sf Tens}}

\def\cot{_{\<\>}}
\def\ten{_{\tensor}}
\def\+{^{\dagger}}
\def\etal{~{\it et~al.}}

\def\Cocont{{\sf Cocont}}

\def\Rel{{\sf Rel}}

\def\:{\colon}
\def\1{{\bf 1}}

\def\2{{\bf 2}}
\def\3{{\bf 3}}

\def\QUANT{{\sf QUANT}}

\def\skel{_{\sf skel}}

\def\op{^{\sf op}}

\def\dom{{\sf dom}}
\def\cod{{\sf cod}}

\def\Sup{{\sf Sup}}
\def\Cat{{\sf Cat}}

\def\Dist{{\sf Dist}}

\def\Cat{{\sf Cat}}

\def\Map{{\sf Map}}

\def\Q{{\cal Q}}
\def\A{{\cal A}}
\def\B{{\cal B}}

\def\P{{\cal P}}
\def\V{{\cal V}}
\def\W{{\cal W}}

\def\F{{\cal F}}
\def\colim{\mathop{\rm colim}}
\def\lim{\mathop{\rm lim}}
\def\bbA{\mathbb{A}}
\def\bbB{\mathbb{B}}
\def\bbC{\mathbb{C}}

\def\bbI{\mathbb{I}}

\def\tensor{\otimes}
\def\<{\langle}
\def\>{\rangle}

\title{Categorical structures enriched in a quantaloid: \\ tensored and cotensored categories}
\author{Isar Stubbe\footnote{D\'epartement de Math\'ematique, Universit\'e de Louvain, Chemin du Cyclotron 2,
1348 Louvain-la-Neuve (Belgique), {\tt
i.stubbe@math.ucl.ac.be}.}}
\date{May 10, 2004\footnote{Corrected some typos on September 18, 2004.}}

\begin{document}

\maketitle

\begin{abstract} Our subject is that of categories, functors and distributors enriched in a base quantaloid $\Q$. We show how cocomplete $\Q$-categories are precisely those which are tensored and conically cocomplete, or alternatively, those which are tensored, cotensored and order-cocomplete. Bearing this in mind, we analyze how $\Sup$-valued homomorphisms on $\Q$ are related to $\Q$-categories. With an appendix on action, representation and variation.
\end{abstract}

\section{Introduction}
The definition of ``category enriched in a bicategory $\W$'' is as old as the definition of bicategory itself [B\'enabou, 1967]. Taking a $\W$ with only one object gives a monoidal category; for symmetric closed monoidal $\V$ the theory of $\V$-categories is well known [Kelly, 1982].
But also categories enriched in a $\W$ with more than one object are interesting. [Walters, 1981] observed that sheaves on a locale give rise to bicategory-enriched categories: ``variation'' (sheaves on a locale $\Omega$) is related to ``enrichment'' (categories enriched in $\Rel(\Omega)$). This insight was further developed in [Walters, 1982] and [Betti\etal, 1983]. Later [Gordon and Power, 1997, 1999] complemented this work, stressing the important r\^ole of tensors in bicategory-enriched categories. 
\par
Here we wish to discuss ``variation and enrichement'' in the case of a base quantaloid (a $\Sup$-enriched category). This is, of course, a particular case of the above, but we believe that it is also of particular interest; many examples of bicategory-enriched categories (like Walters') are really quantaloid-enriched. Since in a quantaloid $\Q$ every diagram of 2-cells commutes, many coherence issues disappear, so the theory of $\Q$-enriched categorical structures is very transparent. Moreover, by definition a quantaloid $\Q$ has stable local colimits, hence (by local smallness) it is closed; this is of great help (to say the least) when working with $\Q$-categories. The theory of quantaloids is documented in [Rosenthal, 1996], and [Stubbe, 2004] provides a reference for all the necessary definitions and basic facts from $\Q$-category theory that will be needed further on.
\par
Our starting point here is the notion of weighted colimit in a $\Q$-category $\bbC$ [Kelly, 1982; Street, 1983]. Two particular cases of such weighted colimits are tensors and conical colimits; then $\bbC$ is cocomplete (i.e.~it admits all weighted colimits) if and only if it is tensored and has all conical colimits [Kelly, 1982; Gordon and Power, 1999] (see also \ref{x3} below). But we may consider the family of ordered sets of objects of the same type in $\bbC$; we call $\bbC$ ``order-cocomplete'' when these ordered sets admit arbitrary suprema. This is a weaker requirement than for $\bbC$ to have conical colimits, but for cotensored $\bbC$ they coincide.
Now $\bbC$ is cocomplete if and only if it is tensored, cotensored and order-cocomplete (as in \ref{x8.1}). Put differently, for a tensored and cotensored $\Q$-category $\bbC$, order-theoretical content (suprema) can be ``lifted'' to $\Q$-categorical content (weighted colimits). 
\par
Then a section is devoted to adjunctions. We see how, at least for tensored $\Q$-categories, order-adjunctions can be ``lifted'' to $\Q$-enriched adjunctions, and how (co)tensoredness may be characterized by enriched adjunctions (analogously to $\V$-categories). As a result, for a tensored $\bbC$, its cotensoredness is equivalent to certain order-adjunctions (cf.~\ref{x16.1}). 
\par
With this in mind we analyze the basic biequivalence between tensored $\Q$-enriched categories and closed pseudofunctors on $\Q\op$ with values in $\Cat(\2)$ (as in \ref{x20}, a particular case of results in [Gordon and Power, 1997]). A finetuned version thereof (in \ref{x37}) says that right $\Q$-modules are the same thing as cocomplete $\Q$-enriched categories.

\section{More on weighted (co)limits}

Throughout $\Q$ denotes a small quantaloid, and our $\Q$-categories have a small set of objects. All notations are as in [Stubbe, 2004].

\subsection*{(Co)tensors}

Let $\bbC$ be a $\Q$-category. For a $\Q$-arrow $f\:X\to Y$ and an object $y\in\bbC_0$ of type $ty=\cod(f)=Y$, the {\em tensor} of $y$ and $f$ is by definition the $(f)$-weighted colimit of $\Delta y$; it will be denoted $y\tensor f$. Thus, whenever it exists, $y\tensor f$ is the (necessarily essentially unique) object of $\bbC$ (necessarily of type $t(y\tensor f)=\dom(f)$) such that 
$$\mbox{for all $z\in\bbC$, }\bbC(y\tensor f,z)=\Big[f,\bbC(y,z)\Big]\mbox{ in }\Q.$$
\par
A {\em cotensor} in $\bbC$ is a tensor in the $\Q\op$-category $\bbC\op$; in elementary terms, for an arrow $f\:X\to Y$ in $\Q$ and an object $x\in\bbC$ of type $tx=\dom(f)=X$, the cotensor of $f$ and $x$, denoted $\<f,x\>$, is -- whenever it exists -- the object of $\bbC$ of type $t\<f,x\>=\cod(f)$ with the universal property that
$$\mbox{for all $z\in\bbC$, }\bbC(z,\<f,x\>)=\Big\{f,\bbC(z,x)\Big\}\mbox{ in }\Q.$$
Thus, $\<f,x\>$ is the $(f)$-weighted limit of $\Delta x$.
\par
A $\Q$-category $\bbC$ is {\em tensored} when for all
$f\in \Q$ and $y\in\bbC_0$ with $ty=\cod(f)$, the tensor $y\tensor f$ exists; and $\bbC$ is {\em cotensored} when $\bbC\op$ is tensored. 
\par
When making a theory of (small) tensored $\Q$-categories, there are some size issues to address, as the following indicates.
\begin{lemma}\label{x0}
A tensored $\Q$-category has either no objects at all, or at least one object of type $X$ for each $\Q$-object $X$. 
\end{lemma}
\proof
The empty $\Q$-category is trivially tensored. Suppose that $\bbC$ is non-empty and tensored; say that there is an object $y$ of type $ty=Y$ in $\bbC$. Then, for any $\Q$-object $X$ the tensor of $y$ with the zero-morphism $0_{X,Y}\in\Q(X,Y)$ must exist, and is an object of type $X$ in $\bbC$.
\endofproof
This motivates once more why {\em we work over a small base quantaloid $\Q$}.
\begin{example}
The two-element Boolean algebra is denoted $\2$; we may view it as a one-object quantaloid so that $\2$-categories are ordered sets, functors are order-preserving maps, and distributors are ideal relations.
A non-empty $\2$-category, i.e.~a non-empty order, is tensored if and only if it has a bottom element, and cotensored if and only if it has a top element. \end{example}
\begin{example}
For any object $Y$ in a quantaloid $\Q$, $\P Y$ denotes the $\Q$-category of contravariant presheaves on the one-object $\Q$-category $*_Y$ whose hom-arrow is $1_Y$. It is cocomplete, thus complete, thus both tensored and cotensored. For an object $f\in\P Y$ of type $tf=X$ (i.e.~a $\Q$-arrow $f\:X\to Y$) and a $\Q$-arrow $g\:U\to X$, $f\tensor g=f\circ g\:U\to Z$ seen as object of type $U$ in $\P Y$. For $h\:X\to V$, $\<h,f\>=\{h,f\}\:V\to Y$, an object of type $V$ in $\P Y$. Similarly, $\P\+ X$ is the $\Q$-category of covariant presheaves on $*_X$; for $f\:X\to Y$, $k\:Y\to M$ and $l\:N\to Y$, $f\tensor l=[l,f]$ and $\<k,f\>=k\circ f$ in $\P\+ X$.
\end{example}

\subsection*{Conical (co)limits}

A $\Q$-category $\bbC$ has an underlying order $(\bbC_0,\leq)$: put $x'\leq x$ whenever both these objects are of the same type, say $tx=tx'=X$, and $1_X\leq\bbC(x',x)$. Conversely, on an ordered set $(A,\leq)$ we may
consider the free $\Q(X,X)$-category $\bbA$:
\begin{itemize}
\item $\bbA_0=A$, all objects are of type $X$;
\item $\bbA(a',a)=\left\{\begin{array}{cl}
1_X & \mbox{if }a'\leq a,\\
0_{X,X} & \mbox{otherwise.}
\end{array}\right.$
\end{itemize}
To give a functor $F\:\bbA\to\bbC$ is to give objects $Fa$, $Fa'$, ... of type $X$ in $\bbC$, such that $Fa'\leq Fa$ in the underlying order of $\bbC$ whenever $a'\leq a$ in $(A,\leq)$.
Consider furthermore the weight $\phi\:*_X\dist\bbA$
whose elements are $\phi(a)=1_X$ for all $a\in\bbA_0$. 
The $\phi$-weighted colimit of $F\:\bbA\to\bbC$ (which may or may not exist) is the {\em conical colimit of $F$}. (Notwithstanding the adjective ``conical'', this is still a weighted colimit!)
A {\em conically cocomplete} $\Q$-category is one that admits all conical colimits\footnote{Analogously to \ref{x0}, a conically cocomplete $\Q$-category $\bbC$ has, for each $\Q$-object $X$, at least one object of type $X$. Indeed,
the conical colimit on the empty functor from the empty free $\Q(X,X)$-category into $\bbC$ is an object of type $X$ in $\bbC$.}. 
\par
The dual notions are those of {\em conical limit} and {\em conically complete $\Q$-category}. We do not bother spelling them out.
\par
The following will help us calculate conical colimits.
\begin{proposition}\label{x1}
Consider a free $\Q(X,X)$-category $\bbA$ and a functor $F\:\bbA\to\bbC$.
An object $c\in\bbC_0$, necessarily of type $tc=X$, is the conical colimit of $F$ if and only if 
$\bbC(c,-)=\bigwedge_{a\in\bbA_0}\bbC(Fa,-)$ in $\Dist(\Q)(\bbC,*_X)$.
\end{proposition}
\proof
For the conical colimit weight $\phi\:*_X\dist\bbA$, 
$\phi(a)=1_X$ for all $a\in\bbA$, thus $c=\colim(\phi,F)$ if and only if
\begin{eqnarray*}
\bbC(c,-)
 & = & \Big[\phi,\bbC(F-,-)\Big] \\
 & = & \bigwedge_{a\in\bbA_0}\Big[\phi(a),\bbC(Fa,-)\Big] \\
 & = & \bigwedge_{a\in\bbA_0}\Big[1_X,\bbC(Fa,-)\Big] \\
 & = & \bigwedge_{a\in\bbA_0}\bbC(Fa,-).
\end{eqnarray*}
\endofproof
In the proof above, to pass from the first line to the second in the series of equations, we used the explicit formula for liftings in the quantaloid $\Dist(\Q)$: in general, for distributors
$\Theta\:\bbA\dist\bbC$ and
$\Psi\:\bbB\dist\bbC$ between $\Q$-categories, $[\Psi,\Theta]\:\bbA\dist\bbB$ has elements, for $a\in\bbA_0$ and $b\in\bbC_0$,
$[\Psi,\Theta](b,a)=\bigwedge_{c\in\bbC_0}[\Psi(c,b),\Theta(c,a)]$,
where the liftings on the right are calculated in $\Q$.
\par
\begin{proposition}\label{x2}
A $\Q$-category $\bbC$ is conically cocomplete if and only if for any family $(c_i)_{i\in I}$ of objects of $\bbC$, all of the same type, say $tc_i=X$, there exists an object $c$ in $\bbC$, necessarily also of that type, such that $\bbC(c,-)=\bigwedge_{i\in I}\bbC(c_i,-)$ in $\Dist(\Q)(\bbC,*_X)$.
\end{proposition}
\proof
One direction is a direct consequence of \ref{x1}. For the other, given a family $(c_i)_{i\in I}$ of objects of $\bbC$, all of type $tc_i=X$, consider the free $\Q(X,X)$-category $\bbI$ on the ordered set $(I,\leq)$ with 
$i\leq j\iff c_i\leq c_j$ in $\bbC$.
The conical colimit of the functor $F\:\bbI\to\bbC\:i\mapsto c_i$ is an object $c\in\bbC_0$ such that $\bbC(c,-)=\bigwedge_{i\in I}\bbC(c_i,-)$, precisely what we wanted.
\endofproof
In what follows we will often speak of ``the conical (co)limit of a family of objects with the same type'', referring to the construction as in the proof above.
\begin{theorem}\label{x3}
A $\Q$-category $\bbC$ is cocomplete if and only if it is tensored and conically cocomplete.
\end{theorem}
\proof
For the non-trivial implication, the alternative description of conical cocompleteness in \ref{x2} is useful. If $\phi\:*_X\dist\bbC$ is {\em any} presheaf on $\bbC$, then the conical colimit of the family 
$(x\tensor\phi(x))_{x\in\bbC_0}$
is the $\phi$-weighted colimit of $1_{\bbC}$: for this is an object $c\in\bbC_0$ such that
\begin{eqnarray*}
\bbC(c,-)
 & = & \bigwedge_{x\in\bbC_0}\bbC(x\tensor\phi(x),-) \\
 & = & \bigwedge_{x\in\bbC_0}\Big[\phi(x),\bbC(x,-)\Big] \\
 & = & \Big[\phi,\bbC(1_{\bbC}-,-)\Big].
\end{eqnarray*}
Hence $\bbC$ is cocomplete (indeed, it suffices that $\bbC$ admit presheaf-weighted colimits of $1_{\bbC}$).
\endofproof
Tensors and conical colimits allow for a very explicit description of colimits in a cocomplete category.
\begin{corollary}\label{x4}
If $\bbC$ is a cocomplete $\Q$-category, then the colimit of
$$\xymatrix@=15mm{\bbA\ar[r]|{\distsign}^{\Phi}&\bbB\ar[r]^F&\bbC}$$
is the functor $\colim(\Phi,F)\:\bbA\to\bbC$ sending an object $a\in\bbA_0$ to the conical colimit of the family $(Fb\tensor\Phi(b,a))_{b\in\bbB_0}$.
A functor $F\:\bbC\to\bbC'$ between cocomplete $\Q$-categories is cocontinuous if and only if it preserves tensors and conical colimits.
\end{corollary}
In \ref{x8.2} we will discuss a more user-friendly version of the above: we can indeed avoid the {\em conical colimits}, and replace them by suitable {\em suprema}.

\subsection*{A third kind of (co)limit}

It makes no sense to ask for the underlying order $(\bbC_0,\leq)$ of a $\Q$-category $\bbC$ to admit arbitrary suprema: two objects of different type cannot even have an upper bound! So let us now denote $\bbC_X$ for the {\em ordered set of $\bbC$-objects with type $X$} (which is thus the empty set when $\bbC$ has no such objects); in these orders it does make sense to talk about suprema. We will say that $\bbC$ is {\em order-cocomplete} when each $\bbC_X$ admits all suprema\footnote{An order-cocomplete $\Q$-category $\bbC$ has, for each $\Q$-object $X$, at least one object of type $X$. Namely, each $\bbC_X$ contains the empty supremum, i.e.~has a bottom element. So (small) order-cocomplete $\Q$-categories can only exist over a small base quantaloid.}.
\par
The dual notion is that of {\em order-complete $\Q$-category}; but of course ``order-complete'' and ``order-cocomplete'' are always equivalent since each order $\bbC_X$ is {\em small}. Nevertheless we will pedantically use both terms, to indicate whether we take suprema or infima as primitive structure.
\par
\begin{proposition}\label{x5}
Let $\bbC$ be a $\Q$-category. 
The conical colimit of a family $(c_i)_{i\in I}\in\bbC_X$ is also its supremum in $\bbC_X$.
\end{proposition}
\proof
Use that $\bbC(c,-)=\bigwedge\bbC(c_i,-)$ in $\Dist(\Q)(\bbC,*_X)$ for the conical colimit $c\in\bbC_0$ of the given family to see that $c=\bigvee_i c_i$ in $\bbC_X$.
\endofproof
So if $\bbC$ is a conically cocomplete $\Q$-category, then it is also order-cocomplete. The converse is not true in general without extra assumptions.
\begin{example}\label{x5.1}
Consider the $\Q$-category $\bbC$ that has, for each $\Q$-object $X$, precisely one object of type $X$; denote this object as $0_X$. The hom-arrows in $\bbC$ are defined as $\bbC(0_X,0_X)=1_X$ (the identity arrow in $\Q(X,X)$) and $\bbC(0_Y,0_X)=0_{X,Y}$ (the bottom element in $\Q(X,Y)$). Then each $\bbC_X=\{0_X\}$ is a sup-lattice, so $\bbC$ is order-cocomplete. However the conical colimit of the {\em empty} family of objects of type $X$ does not exist as soon as the identity arrows in $\Q$ are not the top elements, or as soon as $\Q$ has more than one object.
\end{example}
\begin{proposition}\label{x6}
Let $\bbC$ be a cotensored $\Q$-category.
The supremum of a family $(c_i)_{i\in I}\in\bbC_X$ is also its conical colimit in $\bbC$.
\end{proposition}
\proof
By hypothesis the supremum $\bigvee_ic_i$ in $\bbC_X$ exists, and by \ref{x5} it is the only candidate to be the wanted conical colimit. Thus we must show that $\bbC(\bigvee_ic_i,-)=\bigwedge_i\bbC(c_i,-)$. But this follows from the following adjunctions between orders:
$$\mbox{for any $y\in\bbC_Y$, }\bbC_X\xymatrix@=20mm{ \ar@{}[r]|{\perp}\ar@<1mm>@/^3mm/[r]^{\bbC(-,y)} & \ar@<1mm>@/^3mm/[l]^{\<-,y\>}}\Q(Y,X)\op\mbox{ in }\Cat(\2).
$$
A direct proof\footnote{Actually these adjunctions in $\Cat(\2)$ follow from adjunctions in $\Cat(\Q)$ which are due to the cotensoredness of $\bbC$---see~\ref{x9.0}.} for this adjunction is easy: one uses cotensors in $\bbC$ to see that, for any $x\in\bbC_X$,
\begin{itemize}
\item $1_X\leq\Big\{\bbC(x,y),\bbC(x,y)\Big\}=\bbC(x,\<\bbC(x,y),y\>)$ hence $x\leq\<\bbC(x,y),y\>$ in $\bbC_X$;
\item $1_X\leq\bbC(\<f,y\>,\<f,y\>)=\Big\{f,\bbC(\<f,y\>,y)\Big\}$ hence $\bbC(\<f,y\>,y)\leq\op f$ in $\Q(Y,X)$.
\end{itemize}
Any left adjoint between orders preserves all suprema that happen to exist, so for any $y\in\bbC_Y$, $\bbC(\bigvee_ic_i,y)=\bigwedge_i\bbC(c_i,y)$ in $\Q(Y,X)$, hence -- since infima of distributors are calculated elementwise -- $\bbC(\bigvee_ic_i,-)=\bigwedge_i\bbC(c_i,-)$ in $\Dist(\Q)(\bbC,*_X)$.
\endofproof
So if $\bbC$ is cotensored and order-cocomplete, then it is also conically cocomplete.
Put differently, a cotensored $\Q$-category is conically cocomplete if and only if it is order-cocomplete.
Dually, a tensored category is conically complete if and only if it is order-complete. So... 
\begin{theorem}\label{x8.1}
For a tensored and cotensored $\Q$-category, all notions of completeness and cocompleteness coincide.
\end{theorem}
As usual, for orders the situation is much simpler than for general $\Q$-categories.
\begin{example}
For any $\2$-category (be it {\em a priori} tensored and cotensored or not) all notions of completeness and cocompleteness coincide: an order is order-cocomplete if and only if it is order-complete, but it is then non-empty and has bottom and top element, thus it is tensored and cotensored, thus it is also conically complete and cocomplete, thus also complete and cocomplete {\em tout court}.
\end{example}
\par
In \ref{x4} arbitrary colimits in a cocomplete $\Q$-category are reduced to tensors and conical colimits. But a cocomplete $\Q$-category is always complete too; so in particular cotensored. By cotensoredness the conical colimits may be further reduced to suprema.
\begin{corollary}\label{x8.2}
If $\bbC$ is a cocomplete $\Q$-category, then the colimit of the diagram
$$\xymatrix@=15mm{\bbA\ar[r]|{\distsign}^{\Phi}&\bbB\ar[r]^F&\bbC}$$
is the functor $\colim(\Phi,F)\:\bbA\to\bbC$ sending an object $a\in\bbA_0$ to the supremum of the family $(Fb\tensor\Phi(b,a))_{b\in\bbB_0}$. And a functor $F\:\bbC\to\bbC'$ between cocomplete $\Q$-categories is cocontinuous if and only it preserves tensors and suprema in each of the $\bbC_X$.
\end{corollary}

\section{(Co)tensors and adjunctions}

\subsection*{Adjunctions and adjunctions are two}

An adjunction of functors between $\Q$-categories, like
$$\bbA\xymatrix@=20mm{ \ar@{}[r]|{\perp}\ar@<1mm>@/^3mm/[r]^F & \ar@<1mm>@/^3mm/[l]^G}\bbB,$$
means that $G\circ F\geq 1_{\bbA}$ and $F\circ G\leq 1_{\bbB}$ in $\Cat(\Q)$. Since functors are type-preserving, this trivially implies adjunctions
$$\mbox{for any $\Q$-object $X$, }\bbA_X\xymatrix@=20mm{ \ar@{}[r]|{\perp}\ar@<1mm>@/^3mm/[r]^{F} & \ar@<1mm>@/^3mm/[l]^{G}}\bbB_X\mbox{ in }\Cat(\2).$$
Now we are interested in the converse: how do adjunctions in $\Cat(\2)$ determine adjunctions in $\Cat(\Q)$? The pertinent result is the following.
\begin{theorem}\label{x12}
Let $F\:\bbA\to\bbB$ be a functor between $\Q$-categories, with $\bbA$ tensored. Then the following are equivalent:
\begin{enumerate}
\item $F$ is a left adjoint in $\Cat(\Q)$;
\item $F$ preserves tensors and, for all $\Q$-objects $X$, $F\:\bbA_X\to\bbB_X$ is a left adjoint in $\Cat(\2)$.
\end{enumerate}
\end{theorem}
\proof
One direction is trivial. For the other, denote the assumed adjunctions in $\Cat(\2)$ as
$$\xymatrix@=20mm{\bbA_X\ar@{}[r]|{\perp}\ar@<1mm>@/^3mm/[r]^{F} & \bbB_X\ar@<1mm>@/^3mm/[l]^{G_X}}\mbox{, one for each $\Q$-object $X$.}$$
First, for any $a\in\bbA_X$ and $b\in\bbB_Y$,
\begin{eqnarray*}
\bbA(a,G_Yb)
 & \leq & \bbB(Fa,FG_Yb) \\
 & = & \bbB(Fa, FG_Yb)\circ 1_Y \\
 & \leq & \bbB(Fa,FG_Yb)\circ\bbB(FG_Yb,b) \\
 & \leq & \bbB(Fa,b).
\end{eqnarray*}
The first inequality holds by functoriality of $F$; to pass from the second to the third line, use the pertinent adjunction $F\dashv G_Y$: $FG_Yb\leq b$ in $\bbB_Y$, so $1_Y\leq\bbB(FG_Yb,b)$. For the converse inequality, use tensors in $\bbA$ and the fact that $F$ preserves them: for $a\in\bbA_X$ and $b\in\bbB_Y$,
\begin{eqnarray*}
\bbB(Fa,b)\leq\bbA(a,G_Yb)
 & \iff & 1_Y\leq\Big[\bbB(Fa,b),\bbA(a,G_Yb)\Big] \\
 & \iff & 1_Y\leq\bbA\Big(a\tensor\bbB(Fa,b),G_Yb\Big) \\
 & \iff & 1_Y\leq\bbB\Big(F(a\tensor\bbB(Fa,b)),b\Big) \\
 & \iff & \bbB\Big(Fa\tensor\bbB(Fa,b),b\Big) \\
 & \iff & \Big[\bbB(Fa,b),\bbB(Fa,b)\Big]
\end{eqnarray*}
which is true. It remains to prove that $G\:\bbB\to\bbA\:b\mapsto G_{tb}b$ is a functor; but for $b\in\bbB_Y$ and $b'\in\bbB_{Y'}$,
\begin{eqnarray*}
\bbB(b',b)
 & = & 1_{Y'}\circ\bbB(b',b) \\
 & \leq & \bbB(FG_{Y'}b',b')\circ\bbB(b',b) \\
 & \leq & \bbB(FG_{Y'}b',b) \\
 & = & \bbA(G_{Y'}b',G_Yb).
\end{eqnarray*}
Here we use once more the suitable $F\dashv G_{Y'}$, but also the composition in $\bbB$ and the equality $\bbB(Fa,b)=\bbA(a,G_{Y}b)$.
\endofproof
\par
In a way, \ref{x12} ressembles \ref{x6}: in both cases $\2$-categorical content is ``lifted'' to $\Q$-categorical content (suprema are ``lifted'' to conical colimits, adjunctions between orders are ``lifted'' to adjunctions between categories), and in both cases the price to pay has to do with (existence and preservation of) (co)tensors.
\par
There is a ``weaker'' version of \ref{x12}: given two functors $F\:\bbA\to\bbB$ and $G\:\bbB\to\bbA$, $F\dashv G$ in $\Cat(\Q)$ if and only if, for each $\Q$-object $X$, $F_X\dashv G_X$ in $\Cat(\2)$. Here one needn't ask $\bbA$ to be tensored nor $F$ to preserve tensors (although it does {\em a posteriori} for it is a left adjoint). But the point is that for this ``weaker'' proposition one {\em assumes the existence} of some functor $G$ and one proves that it is the right adjoint to $F$, whereas in \ref{x12} one {\em proves the existence} of the right adjoint to $F$.
\par
Were we to prove \ref{x12} under the hypothesis that $\bbA$, $\bbB$ are cocomplete $\Q$-categories, we simply could have applied \ref{x8.2}: for such categories, $F\:\bbA\to\bbB$ is left adjoint if and only if it is cocontinuous, if and only if preserves tensors and each $\bbA_X\to\bbB_X\:a\mapsto Fa$ preserves suprema, if and only if it preserves tensors and each $\bbA_X\to\bbB_X\:a\mapsto Fa$ is left adjoint in $\Cat(\2)$ (for each $\bbA_X$ is a cocomplete order). The merit of \ref{x12} is thus to have generalized \ref{x8.2} to the case of a tensored $\bbA$ and an arbitrary $\bbB$.

\subsection*{Adjunctions from (co)tensors, and {\em vice versa}}

\begin{proposition}\label{x9.0}
For a $\Q$-category $\bbC$ and an object $x\in\bbC_X$, all cotensors with $x$ exist if and only if the functor\footnote{In principle, $\bbC(-,x)\:*_X\dist\bbC$ is a covariant presheaf on $\bbC$, i.e.~a distributor; but these correspond precisely to functors from $\bbC$ to the completion of $*_X$, which we denote as $\P\+ X$. We do not notationally distinguish between distributor and functor here.} $\bbC(-,x)\:\bbC\to\P\+ X$ is a left adjoint in $\Cat(\Q)$. In this case its right adjoint is $\<-,x\>\:\P\+ X\to\bbC$.
\end{proposition}
\proof
If for any $f\:X\to Y$ in $\Q$ the cotensor $\<f,x\>$ exists, then $\<-,x\>\:\P\+ X\to\bbC$ is a functor: for $f\:X\to Y$, $f'\:X\to Y'$, i.e.~two objects of $\P\+ X$, 
\begin{eqnarray*}
\P\+ X(f',f)\leq\bbC(\<f',x\>,\<f,x\>) 
 & \iff & \Big\{f,f'\Big\}\leq\Big\{f,\bbC(\<f',x\>,x)\Big\} \\
 & \Longleftarrow & f'\leq\bbC(\<f',x\>,x) \\
 & \iff & 1_{Y'}\leq\bbC(\<f',x\>,\<f',x\>)
\end{eqnarray*}
which is true. And $\bbC(-,x)\dashv\<-,x\>$ holds by the universal property of the cotensor itself.
\par
Conversely, suppose that $\bbC(-,x)\:\bbC\to\P\+ X$ is a left adjoint; let $R_x\:\P\+ X\to\bbC$ denote its right adjoint. Then in particular for all $f\:X\to Y$ in $\Q$, $R_x(f)$ is an object of type $Y$ in $\bbC$, satisfying 
$$\mbox{for all $y\in\bbC$, }\bbC(y,R_x(f))=\P\+ X\Big(\bbC(y,x),f\Big)
=\Big\{f,\bbC(y,x)\Big\},$$
which says precisely that $R_x(f)$ is the cotensor of $x$ with $f$.
\endofproof
In the situation of \ref{x9.0} it follows that
\begin{eqnarray}\label{x10}
 & \mbox{for each $\Q$-object $Z$, }\bbC_Z\xymatrix@=20mm{ \ar@{}[r]|{\perp}\ar@<1mm>@/^3mm/[r]^{\bbC(-,x)} & \ar@<1mm>@/^3mm/[l]^{\<-,x\>}}\Q(X,Z)\op\mbox{ in }\Cat(\2), & \\
 & \mbox{for each $z\in\bbC_Z$, }\bbC(z,x)=\bigwedge\{f\:X\to Z\mbox{ in }\Q\mid z\leq\<f,x\>\mbox{ in }\bbC_Z\}. & 
\end{eqnarray}
\par
The dual version of the above will be useful too: it says that tensors with $y\in\bbC_Y$ exist if and only if $\bbC(y,-)\:\bbC\to\P Y$ is a right adjoint in $\Cat(\Q)$, in which case its left adjoint is $y\tensor-\:\P Y\to\bbC$. 
And then moreover
\begin{eqnarray}
 & \mbox{for each $\Q$-object $Z$, }\bbC_Z\xymatrix@=20mm{ \ar@{}[r]|{\perp}\ar@<-1mm>@/_3mm/[r]_{\bbC(y,-)} & \ar@<-1mm>@/_3mm/[l]_{y\tensor-}}\Q(Z,Y)\mbox{ in }\Cat(\2), & \\
 \label{x10.1}
 & \mbox{for each $z\in\bbC_Z$, }\bbC(y,z)=\bigvee\{f\:Z\to Y\mbox{ in }\Q\mid y\tensor f\leq z\mbox{ in }\bbC_Z\}. &
\end{eqnarray}
\par
Here is a usefull application of the previous results. For any $\Q$-category $\bbC$ the Yoneda embedding $Y\+_{\bbC}\:\bbC\to\P\+\bbC\:c\mapsto\bbC(c,-)$ is a cocontinuous functor; in particular, for any $x\in\bbC_X$ the functor $\bbC(-,x)\:\bbC\to\P\+ X$ preserves tensors. (A direct proof of this latter fact is easy too: for $f\:Y\to Z$ in $\Q$ and $z\in\bbC_Z$, suppose that $z\tensor f$ exists in $\bbC$. Then $\bbC(z\tensor f,x)=[f,\bbC(z,x)]=\bbC(z,x)\tensor f$ in $\P\+ X$, because this is how tensors are calculated in $\P\+ X$.)
\begin{corollary}\label{x14}
If $\bbC$ is a tensored $\Q$-category, then the following are equivalent:
\begin{enumerate}
\item for all $\Q$-objects $X$ and $Y$ and each $x\in\bbC_X$, $\bbC(-,x)\:\bbC_Y\to\Q(X,Y)\op$ is a left adjoint in $\Cat(\2)$;
\item for each $x\in\bbC_X$, $\bbC(-,x)\:\bbC\to\P\+ X$ is a left adjoint in $\Cat(\Q)$;
\item  $\bbC$ is cotensored.
\end{enumerate}
\end{corollary}
\par
In \ref{x9.0} we have results about ``(co)tensoring with a fixed object''; now we are interested in studying ``tensoring with a fixed arrow''. Recall that a tensor is a colimit of which such an arrow is the weight. So we may apply general lemmas on weighted colimits to obtain the following particular results.
\begin{proposition}\label{x15}
Let $\bbC$ denote a $\Q$-category.
\begin{enumerate}
\item For all $y\in\bbC_Y$, $y\tensor 1_Y\cong y$.
\item For $g\:W\to X$ and $f\:X\to Y$ in $\Q$ and $y\in\bbC_Y$, if all tensors involved exist then $y\tensor(f\circ g)\cong(y\tensor f)\tensor g$.
\item for $(f_i\:X\to Y)_{i\in I}$ in $\Q$ and $y\in\bbC_Y$, if all tensors involved exist then $y\tensor(\bigvee_if_i)\cong\bigvee_i(y\tensor f_i)$.
\item For $f\:X\to Y$ in $\Q$ and $y,y'\in\bbC_Y$, if all tensors involved exist then $y\leq y'$ in $\bbC_Y$ implies $y\tensor f\leq y'\tensor f$ in $\bbC_X$.
\end{enumerate}
\end{proposition}
Of course there is a dual version about cotensors, but we do not bother spelling it out. However, there is an interesting interplay between tensors and cotensors.
\begin{proposition}\label{x16}
Let $f\:X\to Y$ be a $\Q$-arrow and suppose that all tensors and all cotensors with $f$ exist in some $\Q$-category $\bbC$. Then
$$\xymatrix@=20mm{\bbC_Y\ar@{}[r]|{\perp}\ar@<1mm>@/^3mm/[r]^{-\tensor f} & \bbC_X\ar@<1mm>@/^3mm/[l]^{\<f,-\>}}\mbox{ in }\Cat(\2).$$
\end{proposition}
\proof
It follows from \ref{x15} (and its dual) that $-\tensor f\:\bbC_Y\to\bbC_X$ and $\<f,-\>\:\bbC_X\to\bbC_Y$ are order-preserving morphisms. Furthermore, for $x\in\bbC_X$ and $y\in\bbC_Y$,
\begin{eqnarray*}
y\tensor f\leq x
 & \iff & 1_X\leq\bbC(y\tensor f,x)=\Big[f,\bbC(y,x)\Big] \\
 & \iff & f\leq\bbC(y,x) \\
 & \iff & 1_Y\leq\Big\{f,\bbC(y,x)\Big\}=\bbC(y,\<f,x\>) \\
 & \iff & y\leq\<f,x\>.
\end{eqnarray*}
\endofproof
We can push this further.
\begin{proposition}\label{x16.1}
A tensored $\Q$-category $\bbC$ is cotensored if and only if, for every $f\:X\to Y$ in $\Q$, $-\tensor f\:\bbC_Y\to\bbC_X$ is a left adjoint in $\Cat(\2)$. In this case, its right adjoint is $\<f,-\>\:\bbC_X\to\bbC_Y$.
\end{proposition}
\proof
Necessity follows from \ref{x16}. As for sufficiency, by \ref{x14} it suffices to show that for all $\Q$-objects $X$ and $Y$ and every $x\in\bbC_X$, 
$$\bbC(x,-)\:\bbC_Y\to\Q(X,Y)\op\:y\mapsto\bbC(x,y)$$
has a right adjoint in $\Cat(\2)$. Denoting, for a $\Q$-arrow $f\:X\to Y$, the right adjoint to $-\tensor f\:\bbC_Y\to\bbC_X$ in $\Cat(\2)$ as $R_f\:\bbC_X\to\bbC_Y$, the obvious candidate right adjoint to $y\mapsto\bbC(x,y)$ is $f\mapsto R_f(x)$. First note that, if $f\leq\op f'$ in $\Q(X,Y)$ then $R_f(x)\tensor f'\leq R_f(x)\tensor f\leq x$ using $-\tensor f\dashv R_f$, which implies by $-\tensor f'\dashv R_{f'}$ that $R_f(x)\leq R_{f'}(x)$: so 
$$R_{(-)}(x)\:\Q(X,Y)\op\to\bbC_Y\:f\mapsto R_f(x)$$ 
preserves order. Further, for $f\in\Q(X,Y)$ and $y\in\bbC_Y$,
\begin{eqnarray*}
\bbC(y,x) \leq\op f 
 & \iff & f\leq \bbC(y,x) \\
 & \iff & y\tensor f\leq x \\
 & \iff & y\leq R_f(x),
\end{eqnarray*}
so indeed $\bbC(x,-)\dashv R_{(-)}(x)$ in $\Cat(\2)$. Now $\bbC$ is tensored and cotensored, so by \ref{x16} it follows that $R_f(x)$ must be $\<f,x\>$ (since both are right adjoint to $-\tensor f$).
\endofproof

\section{Enrichment and variation}

\subsection*{Terminology and notations}

We must introduce some notation. By $\Cat\ten(\Q)$ we denote the full sub-2-category of $\Cat(\Q)$ whose objects are tensored categories, and $\Tens(\Q)$ the sub-2-category whose objects are tensored categories and morphisms are tensor-preserving functors. Similarly we use $\Cat\cot(\Q)$ for the full sub-2-category of $\Cat(\Q)$ whose objects are cotensored categories, and moreover the obvious combination $\Cat_{\tensor,\<\>}(\Q)$. Recall also that $\Cocont(\Q)$ denotes the locally completely ordered 2-category whose objects are cocomplete $\Q$-categories and morphisms are cocontinuous (equivalently, left adjoint) functors; and $\Cocont\skel(\Q)$ denotes its biequivalent full sub-quantaloid whose objects are skeletal.
\begin{example}
$\Cat(\2)$ is the locally ordered 2-category of orders and order preserving maps. $\Cat\ten(\2)$ has orders with bottom element as objects and all order-preserving maps as morphisms, whereas $\Tens(\2)$ has the same objects but the morphisms are required to send bottom onto bottom. $\Cocont(\2)$ is biequivalent to the quantaloid of sup-lattices and sup-morphisms; taking only skeletal $\2$-categories (i.e.~antisymmetric orders) we have $\Cocont\skel(\2)=\Sup$.
\end{example}
\par
Some more notions and notations, now from the realm of ``variation'': Let $\A$ and $\B$ be locally ordered 2-categories (i.e.~$\Cat(\2)$-enriched categories). A {\em pseudofunctor} $\F\:\A\to\B$ is an action on objects and morphisms that respects the local order and such that functoriality holds up to local isomorphism (we needn't require any coherence because our 2-categories are locally ordered). For two such pseudofunctors $\F,\F'\:\A\biar\B$, a {\em lax natural transformation} $\varphi\:\F\tto\F'$ is a family of $\B$-morphisms $(\varphi_X\:\F X\to\F' X)_{X\in\A_0}$ satisfying, for any $f\:X\to Y$ in $\A$, $\F' f\circ\varphi_X\leq\varphi_Y\circ\F f$ in $\B(\F X,\F' Y)$. Such a transformation is {\em pseudonatural} when these inequalities are isomorphisms. Lax natural transformations are ordered componentwise. There are locally ordered 2-categories $\Psd\lax(\A,\B)$, resp.~$\Psd(\A,\B)$, with pseudofunctors as objects and lax natural transformations, resp.~pseudonatural transformations, as arrows.
\par
Now consider a pseudofunctor $\F\:\A\to\Cat(\2)$; it is {\em closed} when, for every $X,Y$ in $\A$ and $x\in\F X$, 
$$\F(-)(x)\:\A(X,Y)\to\F Y\:f\mapsto \F(f)(x)$$
is a left adjoint in $\Cat(\2)$. We write $\ClPsd\lax(\A,\Cat(\2))$ and $\ClPsd(\A,\Cat(\2))$ for the full sub-2-categories of $\Psd\lax(\A,\Cat(\2))$ and $\Psd(\A,\Cat(\2))$ determined by the closed pseudofunctors. 
\par
We will be interested in closed pseudofunctors on the opposite of a quantaloid $\Q$; the closedness of a pseudofunctor $\F\:\Q\op\to\Cat(\2)$ reduces to the fact that, for each $X,Y$ in $\Q$ and $y\in Y$,
\begin{equation}\label{x20.0}
\F(-)(y)\:\Q(X,Y)\to\F X\:y\mapsto\F(f)(y)
\end{equation}
preserves arbitrary suprema (for $\Q(X,Y)$ is a sup-lattice).
When we replace $\Cat(\2)$ by any of its sub-2-categories like $\Cat\ten(\2)$, $\Tens(\2)$ and so on, the closedness condition for  pseudofunctors still makes sense: we will mean precisely that the order-morphisms in (\ref{x20.0}) preserve suprema (i.e.~are left adjoints in $\Cat(\2)$).

\subsection*{The basic biequivalence}

\begin{proposition}\label{x20}
A tensored $\Q$-category $\bbC$ determines a closed pseudofunctor
\begin{equation}\label{x21}
\F_{\bbC}\:\Q\op\to\Cat(\2)\:\Big(f\:X\to Y\Big)\mapsto\Big(-\tensor f\:\bbC_Y\to\bbC_X\Big).
\end{equation}
And a functor $F\:\bbC\to\bbC'$ between tensored $\Q$-categories determines a lax natural transformation 
\begin{equation}\label{x22}
\varphi^F\:\F_{\bbC}\tto\F_{\bbC'}\mbox{ with components }\varphi^F_X\:\bbC_X\to\bbC'_X\:x\mapsto Fx.
\end{equation}
\end{proposition}
\proof
For a tensored $\Q$-category $\bbC$, $\F_{\bbC}$ as in the statement of the proposition is well-defined: each $\bbC_X$ is an order and each $-\tensor f\:\bbC_Y\to\bbC_X$ preserves order (by \ref{x15}). Moreover, this action is pseudofunctorial (again by \ref{x15}). And from (the dual of) \ref{x9.0} we know that, for each $X,Y$ in $\Q$ and $y\in\bbC_Y$, 
$$y\tensor-\:\Q(X,Y)\to\bbC_X\:f\mapsto y\tensor f$$
is a left adjoint; so $\F_{\bbC}$ is a closed pseudofunctor.
\par
A functor $F\:\bbC\to\bbC'$ is a type-preserving mapping $F\:\bbC_0\to\bbC'_0\:x\mapsto Fx$ of objects such that $\bbC(y,x)\leq\bbC'(Fy,Fx)$ for all $x,y\in\bbC_0$. With (\ref{x10.1}), this functor-inequality may be rewritten as
\begin{eqnarray*}
 & & \bbC(y,x)\leq\bbC'(Fy,Fx) \\
 & & \iff\mbox{for any $f\:X\to Y$ in $\Q$, if $y\tensor f\leq x$ in $\bbC_X$ then $Fy\tensor f\leq Fx$ in $\bbC'_X$} \\
 & & \iff\mbox{for any $f\:X\to Y$ in $\Q$, $Fy\tensor f\leq F(y\tensor f)$.}
\end{eqnarray*}
(For the last equivalence, necessity folows by application of the previous sentence to $y\tensor f\leq y\tensor f$, whereas for sufficiency one first notes that $y\tensor f\leq x$ implies anyway that $F(y\tensor f)\leq Fx$ so combined with the assumption this gives $Fy\tensor f\leq Fx$.) Thus, such a functor $F\:\bbC\to\bbC'$ is really just a family of mappings
$\bbC_X\to\bbC'_X\:x\mapsto Fx$, one for each $\Q$-object $X$,
which are all order-preserving (by functoriality of $F$) and 
satisfy furthermore for any $f\:X\to Y$ in $\Q$ and $y\in\bbC_Y$ that $Fy\tensor f\leq F(y\tensor f)$.
Having defined components $\varphi^F_X$ as in (\ref{x22}), this says that $\F_{\bbC'}(f)\circ\varphi^F_Y
\leq \varphi^F_X\circ\F_{\bbC}(f)$, for any 
$f\:X\to Y$ in $\Q$. So $\varphi\:\F_{\bbC}\to\F_{\bbC'}$ is a lax natural transformation.
\endofproof
\begin{theorem}\label{x23}
For any quantaloid $\Q$, the action
\begin{equation}\label{x24}
\Cat\ten(\Q)\to\ClPsd\lax(\Q\op,\Cat(\2))\:
\Big(F\:\bbC\to\bbC'\Big)\mapsto\Big(\varphi^F\:\F_{\bbC}\tto\F_{\bbC'}\Big)
\end{equation}
is an equivalence of 2-categories.
\end{theorem}
\proof
Straightforwardly the action in (\ref{x24}) is functorial: the lax natural transformation corresponding to an identity functor is an idenity lax natural transformation; the lax natural transformation corresponding to the composition of functors is the composition of the lax natural transformations corresponding to each of the functors involved.
\par
Now let $\F\:\Q\op\to\Cat(\2)$ be any closed pseudofunctor; then define a $\Q$-category $\bbC^{\F}$ by:
\begin{itemize}
\item for each $\Q$-object $X$, $\bbC^{\F}_X:=\F X$,
\item for $x\in \bbC^{\F}_X$ and $y\in\bbC^{\F}_Y$, $\bbC^{\F}(y,x)=\bigvee\{f\:X\to Y\mbox{ in }\Q\mid \F(f)(y)\leq x\mbox{ in }\bbC^{\F}_X\}$.
\end{itemize}
The supremum involved is really an expression of the closedness of the pseudofunctor: $x\mapsto\bbC^{\F}(y,x)$ is the right adjoint to $f\mapsto\F(f)(y)$ in $\Cat(\2)$. Then $\bbC^{\F}$ is a tensored $\Q$-category: the tensor of some $f\:X\to Y$ and $y\in \F Y$ is precisely $\F(f)(y)$, by (the dual of) \ref{x9.0}. It is clear that $\F\cong\F_{\bbC^{\F}}$. So far for essential surjectivity of (\ref{x24}).
\par
Finally, given tensored $\Q$-categories $\bbC$ and $\bbC'$, the ordered sets $\Cat\ten(\Q)(\bbC,\bbC')$ and $\Psd\lax(\Q\op,\Cat(\2))(\F_{\bbC},\F_{\bbC'})$ are isomorphic: a functor $F\:\bbC\to\bbC'$ between (tensored) $\Q$-categories is completely determined by its action on objects, hence by the family of (order-preserving) mappings $\bbC_X\to\bbC'_X\:x\mapsto Fx$, hence by the components of the corresponding transformation $\varphi^F\:\F_{\bbC}\tto\F_{\bbC'}$. From the proof of \ref{x20} it is clear that $F$ is a functor if and only if $\varphi^F$ is lax natural (thanks to tensoredness of $\bbC$ and $\bbC'$). Furthermore, to say that $F\leq G\:\bbC\biar\bbC'$ in $\Cat(\Q)$ means that, for any $\Q$-object $X$ and any $x\in\bbC_X$, $Fx\leq Gx$ in $\bbC'_X$.
For the lax natural transformations $\varphi^F,\varphi^G$ corresponding to $F,G$ this is really the same thing as saying that
$\varphi^F_X\leq\varphi^G_X$ in $\Cat(\2)$,
in other words, $\varphi^F\leq\varphi^G$ as arrows between (closed) pseudofunctors.
\endofproof
It follows from \ref{x0} and \ref{x23} that
a closed pseudofunctor $\F\:\Q\op\to\Cat(\2)$ either has all of the $\F X$ empty, or none of them. A direct proof is easy too (it is of course a transcription of \ref{x0} modulo the equivalence in \ref{x23}): if $y\in\F Y$, then $\F(0_{X,Y})(y)\in\F X$, where $0_{X,Y}\in\Q(X,Y)$ is the bottom element. So as soon as one of the $\F X$ is non-empty, all of them are. And the empty pseudofunctor is trivially closed.

\subsection*{Finetuning}

Here are some seemingly innocent specifications concerning the 2-functor in \ref{x23}.
\begin{lemma}\label{x25}
Any closed pseudofunctor $\F\:\Q\op\to\Cat(\2)$ lands in $\Cat\ten(\2)$. And any lax natural transformation $\varphi\:\F\tto\F'\:\Q\to\Cat(\2)$ between closed pseudofunctors has components in $\Cat\ten(\2)$ rather than $\Cat(\2)$.
\end{lemma}
\proof
For any closed pseudofunctor $\F\:\Q\op\to\Cat(\2)$, for every $X$ in $\Q$ and $x\in\F X$, $\F(-)(x)\:\Q(X,X)\to\F X$ preserves all suprema, thus in particular the empty supremum, i.e.~the bottom element $0_{X,X}\in\Q(X,X)$. This implies that every non-empty $\F X$ must have a bottom element. Thus $\F$ lands in $\Cat\ten(\2)$ rather than $\Cat(\2)$.
But precisely because of this, 
the components $\varphi_X\:\F X\to\F' X$ of a lax natural transformation $\varphi^F\:\F_{\bbC}X\to\F_{\bbC'}X$ live in $\Cat\ten(\2)$ rather than $\Cat(\2)$.
\endofproof
From this proof it follows that, for a closed pseudofunctor $\F\:\Q\op\to\Cat\ten(\2)$, the bottom element in a non-empty order $\F X$ may be calculated as: $0_X:=\F(0_{X,X})(x)$, where $x$ is an arbitrary element in $\F X$. This allows for the following.
\begin{lemma}\label{x26.1}
A pseudonatural transformation $\varphi\:\F\tto\F'\:\Q\op\to\Cat\ten(\2)$ between closed pseudofunctors has components in $\Tens(\2)$.
\end{lemma}
\proof
If $\F X$ is non-empty, take any $x\in\F X$, then by pseudonaturality of $\varphi$,
$$\varphi_X(0_X)=\varphi_X(\F(0_{X,X})(x))\cong\F'(0_{X,X})(\varphi_X(x))=0'_X.$$
So each component $\varphi_X\:\F X\to\F'X$, {\em a priori} in $\Cat\ten(\2)$, preserves the bottom element if there is one, thus lives in $\Tens(\2)$.
\endofproof
\begin{lemma}\label{x26}
Any closed pseudofunctor $\F\:\Q\op\to\Map(\Cat\ten(\2))$ actually lands in $\Map(\Cat_{\tensor,\<\>}(\2))$.
\end{lemma}
\proof
Taking an arbitrary $x\in\F X$ (presumed non-emtpy), $\F(0_{X,X})^*(x)$ gives the top element of $\F X$. Here $\F(0_{X,X})^*$ denotes the right adjoint to $\F(0_{X,X})$ in $\Cat\ten(\2)$.
So each $\F X$ is an object of $\Cat_{\tensor,\<\>}(\2)$ rather than $\Cat\ten(\2)$.
\endofproof
\par
Now we can apply all this to finetune \ref{x23}.
\begin{proposition}\label{x27}
Let $\bbC$ be a tensored $\Q$-category.
\begin{enumerate}
\item\label{x28.0} The associated pseudofunctor $\F_{\bbC}\:\Q\op\to\Cat(\2)$ factors through $\Cat\ten(\2)$. 
\item\label{x28} $\bbC$ is moreover cotensored if and only if $\F_{\bbC}$ factors through $\Map(\Cat_{\tensor,\<\>}(\2))$.
\item\label{x29} $\bbC$ is cocomplete if and only if $\F_{\bbC}$ factors through $\Cocont(\2)$.
\item\label{x30} $\bbC$ is skeletal and cocomplete if and only if $\F_{\bbC}$ factors through $\Cocont\skel(\Q)$.
\end{enumerate} 
\end{proposition}
\proof
(\ref{x28.0}) Is the content of \ref{x25}.
\par
(\ref{x28}) Is a combination of \ref{x16.1}, \ref{x23} and \ref{x26}.
\par
(\ref{x29}) By \ref{x8.1} a tensored and cotensored $\bbC$ is cocomplete if and only if it is order-cocomplete, i.e.~each $\bbC_X$ is a cocomplete order. Now apply (\ref{x28}), recalling that $\Cocont(\2)$ is precisely the full sub-2-category of $\Map(\Cat_{\tensor,\<\>}(\2))$ determined by the (order-)cocomplete objects.
\par
(\ref{x30}) Is a variation on (\ref{x29}): a $\Q$-category $\bbC$ is skeletal if and only if each $\bbC_X$ is an antisymmetric (i.e.~skeletal) order.
\endofproof
\begin{proposition}\label{x31}
Let $F\:\bbC\to\bbC'$ be a functor between tensored $\Q$-categories.
\begin{enumerate}
\item\label{x32.0} The corresponding lax natural transformation 
$\varphi^F\:\F_{\bbC}\tto\F_{\bbC'}$ has components in $\Cat\ten(\2)$.
\item\label{x32} $F$ is tensor-preserving if and only if $\varphi^F$ is pseudonatural.
\item\label{x33} $F$ is left adjoint if and only if $\varphi^F$ is pseudonatural and its components are in $\Map(\Cat\ten(\2))$.
\item\label{x34} If $\bbC$ and $\bbC'$ are moreover cotensored, then $F$ is left adjoint if and only if $\varphi^F$ is pseudonatural and its components are in $\Map(\Cat_{\tensor,\<\>}(\2))$.
\item\label{x35} If $\bbC$ and $\bbC'$ are cocomplete, then $F$ is left adjoint if and only if $\varphi^F$ is pseudonatural and its components are in $\Cocont(\2)$.
\item\label{x36} If $\bbC$ and $\bbC'$ are skeletal and cocomplete, then $F$ is left adjoint if and only if $\varphi^F$ is pseudonatural and its components are in $\Cocont\skel(\2)$.
\end{enumerate}
\end{proposition}
\proof
(\ref{x32.0}) Is the content of \ref{x25}.
\par
(\ref{x32}) To say that $F\:\bbC\to\bbC'$ preserves tensors, means that
for any $f\:X\to Y$ in $\Q$ and $y\in\bbC_Y$, $F(y\tensor f)\cong Fy\tensor f$ in $\bbC_X$. In terms of the transformation $\varphi^F$ this means that
$\varphi^F_X\circ \F_{\bbC}(f)\cong \F_{\bbC'}\circ\varphi^F_Y$
instead of merely the inequality ``$\geq$''; hence it is pseudonatural instead of merely lax natural.
\par
(\ref{x33}) By \ref{x12} and the previous point.
\par
(\ref{x34}) Is a variation on the previous point, using \ref{x27} (\ref{x28}).
\par
(\ref{x35}) Follows from (\ref{x33}), taking into account that all $\bbC_X$ and $\bbC'_X$ are cocomplete orders.
\par
(\ref{x36}) Is a variation on (\ref{x35}).
\endofproof
We may now state our conclusion.
\begin{theorem}\label{x37}
The equivalence in \ref{x23} reduces to the following equivalences of locally ordered 2-categories:
\begin{enumerate}
\item $\Cat\ten(\Q)\simeq\ClPsd\lax(\Q\op,\Cat\ten(\2))$,
\item $\Tens(\Q)\simeq\ClPsd(\Q\op,\Tens(\2))$,
\item $\Map(\Cat_{\tensor,\<\>}(\Q))
\simeq\ClPsd(\Q\op,\Map(\Cat_{\tensor,\<\>}(\2)))$,
\item $\Cocont(\Q)\simeq\ClPsd(\Q\op,\Cocont(\2))$,
\item $\Cocont\skel(\Q)\simeq\ClPsd(\Q\op,\Cocont\skel(\2))$.
\end{enumerate}
\end{theorem}
Actually, $\Cocont\skel(\2)=\Sup$ and a closed pseudofunctor from $\Q\op$ to $\Sup$ is really a quantaloid homomorphism; moreover, $\Cocont\skel(\Q)$ is biequivalent to $\Cocont(\Q)$. So we may end with the following.
\begin{corollary}\label{x37.1}
The quantaloid of right $\Q$-modules (cf.~\ref{x38}) is biequivalent to the locally cocompletely ordered category of cocomplete $\Q$-categories and cocontinuous functors: $\QUANT(\Q\op,\Sup)\simeq\Cocont(\Q)$.
\end{corollary}

\section{Appendix: action, representation and variation}

Let $K$ denote a 
quantale, i.e.~a one-object quantaloid. Now thinking of $K$ as a monoid in $\Sup$, let ``unit'' and
``multiplication'' in $K$ (the single identity arrow and the composition in the one-object quantaloid) correspond to sup-morphisms
$\varepsilon\:I\to K$ and $\gamma\:K\tensor K\to K$. A {\em right
action} of $K$ on some sup-lattice $M$ is a sup-morphism 
$\phi\:M\tensor K\to M$ such that the diagrams 
$$\xymatrix@=15mm{
M\tensor K\tensor\ar[r]^{1\tensor\gamma}\ar[d]_{\phi\tensor 1_K} K & M\tensor K\ar[d]^{\phi} & M\tensor I\ar[l]_{1\tensor\varepsilon} \\
M\tensor K\ar[r]_{\phi} & M\ar@{=}[ur]}$$
commute (we don't bother writing the associativity and unit isomorphisms in the symmetric monoidal closed category $\Sup$);
$(M,\phi)$ is then said to be a {\em right $K$-module}.
In elementary terms we have a set-mapping $M\times K\to M\:(m,f)\mapsto \phi(m,f)$, preserving suprema in both variables, and such that (with obvious notations) 
$$\phi(m,1)=m\mbox{ and }\phi(m,g\circ f)=\phi(\phi(m,g),f).$$
By closedness of $\Sup$, to the sup-morphism $\phi\:M\tensor K\to M$ corresponds a unique sup-morphism
$\bar{\phi}\:K\to\Sup(M,M)$. In terms of elements, this $\bar\phi$ sends every $f\in K$ to the sup-morphism $\phi(-,f)\:M\to M$; it satisfies 
$$\bar\phi(1)=1_M\mbox{ and }\bar\phi(g\circ f)=\bar\phi(f)\circ\bar\phi(g).$$ That is to say, $\bar\phi\:K\to\Sup(M,M)$ is a {\em reversed representation}
of the quantale $K$ by endomorphisms on the sup-lattice $M$: a
homomorphism of quantales that reverses the multiplication (where $\Sup(M,M)$ is endowed with composition as binary operation and the identity morphism $1_M$ as unit to form a
quantale). Recalling that $K$ is a one-object quantaloid $\Q$, such a multiplication-reversing homomorphism $\bar\phi\:K\to\Sup(M,M)$ is really a $\Sup$-valued
quantaloid homomorphism $F\:\Q\op\to\Sup\:*\mapsto M, f\mapsto
\bar\phi(f)$. 
\par
In the same way it can be seen that morphisms between modules correspond to $\Sup$-enriched natural transformations between $\Sup$-presheaves. Explicitly, for two right modules $(M,\phi)$ and $(N,\psi)$, a module-morphism $\alpha\:M\to N$ is a sup-morphism that makes
$$\xymatrix@=15mm{M\tensor K\ar[d]_{\phi}\ar[r]^{\alpha\tensor 1_K} & N\tensor K\ar[d]^{\psi}\\
M\ar[r]_{\alpha} & N}$$
commute. In elementary terms, such a sup-morphism $\alpha\:M\to N\:m\mapsto \alpha(m)$ satisfies
$$\alpha(\phi(m,f))=\psi(\alpha(m),f).$$
By adjunction -- and with notations as above -- this gives for any $f\in K$ the commutative square
$$\xymatrix@=15mm{M\ar[d]_{\bar\phi(f)}\ar[r]^{\alpha} & N\ar[d]^{\bar\psi(f)} \\
M\ar[r]_{\alpha} & N}$$
which expresses precisely the naturality of $\alpha$ viewed as (single) component of a natural transformation $\alpha\:F\tto G$, where $F,G\:\Q\op\biar\Sup$ denote the homomorphisms corresponding to $M$ and $N$.
\par
Conclusively,
actions, representations and $\Sup$-presheaves are essentially the same
thing. The point now is that the latter presentation straightforwardly makes sense for any quantaloid, and
not just those with only one object.
\begin{definition}\label{x38}
A {\em right $\Q$-module $M$} is a
homomorphism $M\:\Q\op\to\Sup$. And a {\em module-morphism $\alpha\:M\tto N$} between two right $\Q$-modules $M$ and $N$ is an enriched natural transformation between these homomorphisms.
\end{definition}
That is to say, $\QUANT(\Q\op,\Sup)$ is the quantaloid of right $\Q$-modules\footnote{We have chosen here to work with {\em right} actions, {\em reversed} representations, and {\em contravariant} $\Sup$-presheaves. Clearly {\em left} actions correspond to {\em straight} representations and to {\em covariant} $\Sup$-valued presheaves.}.

\end{document}